\documentclass[a4paper,11pt]{amsart}

\usepackage{amssymb}
\usepackage[english]{babel}
\usepackage{amsmath}
\usepackage{amsthm}
\usepackage{verbatim}
\usepackage[all]{xy}
\usepackage{enumerate}
\usepackage{setspace}

\newtheorem{theorem}{Theorem}[section]

\newtheorem{e-proposition}[theorem]{Proposition}
\newtheorem{corollary}[theorem]{Corollary}
\newtheorem{e-definition}[theorem]{Definition\rm}

\newcommand{\bn}{\Bbb N}

\newcommand{\nphi}{\mathfrak{n}_{\varphi}}

\newcommand {\Pcone} {\mathcal{P}}

\newcommand {\mlg} {\mathsf{M}}
\newcommand {\nlg} {\mathsf{N}}

\newcommand{\Hil}{\mathsf{H}}

\newenvironment{rlist}
{

\begin{enumerate}}
{\end{enumerate}}

\newcommand{\la}{\langle}
\newcommand{\ra}{\rangle}

\newcommand{\ot}{\otimes}

\title[Haagerup property for arbitrary von Neumann algebras]{Generalisations of the Haagerup approximation  property to arbitrary von Neumann algebras}

\author{Martijn Caspers}
 \address{M. Caspers, Fachbereich Mathematik und Informatik der Universit\"at M\"unster,
Einsteinstrasse 62,
48149 M\"unster, Germany}
 \email{martijn.caspers@uni-muenster.de}
\thanks{MC is supported by the grant SFB 878 ``{\it Groups, geometry and actions}''}

\author{Rui Okayasu}
\address{R. Okayasu, Department of Mathematics Education, Osaka Kyoiku University, Osaka
582-8582, Japan}
\email{rui@cc.osaka-kyoiku.ac.jp}
\thanks{RO was partially supported by JSPS KAKENHI Grant Number 25800065.}

\author{Adam Skalski}
\address{A. Skalski,  Institute of Mathematics of the Polish Academy of Sciences,
ul.~\'Sniadeckich 8, 00--956 Warszawa, Poland\\ \newline Faculty of Mathematics, Informatics and Mechanics, University of Warsaw, ul.~Banacha 2,
02-097 Warsaw, Poland}
\email{a.skalski@impan.pl}
\thanks{AS was partially supported by the Iuventus Plus grant IP2012 043872.}

\author{Reiji Tomatsu}
\address{R. Tomatsu, Department of Mathematics, Hokkaido University, Hokkaido 060-0810,
JAPAN}
\email{tomatsu@math.sci.hokudai.ac.jp}
\thanks{RT was partially supported by JSPS KAKENHI Grant Number 24740095.}

\begin{document}

\maketitle

\selectlanguage{english}

\begin{abstract}
The notion of the Haagerup approximation property, originally introduced for von Neumann algebras equipped with a faithful normal tracial state, is generalized to arbitrary von Neumann algebras. We discuss two equivalent characterisations, one in terms of the standard form and the other in terms of the approximating maps with respect to a fixed faithful normal semifinite weight. Several stability properties, in particular regarding the crossed product construction are established and certain examples are introduced.

\end{abstract}






\section{Introduction}
\label{sec-intro}
The Haagerup property for a locally compact group has its origins in the celebrated paper \cite{HaaFree}, where U. Haagerup proved that the length function on the free group is conditionally negative definite. This fact, i.e.\ the existence on a group $G$ of a proper continuous conditionally negative definite function can be viewed as a natural relaxation of   amenability of $G$ and has later proved to be a natural and influential notion in geometric group theory, dynamical systems and operator algebras (see \cite{book}). In particular already in 1983 M. Choda showed in \cite{Cho} that if $G$ is discrete, then it has the  Haagerup property if and only if the group von Neumann algebra of $G$ equipped with its natural trace has a certain von Neumann algebraic approximation property, which came to be known as the \emph{Haagerup approximation property} for a von Neumann algebra with a faithful normal tracial state. The latter property was also motivated by the study of cocycle actions on von Neumann algebras (\cite{CoJ}). Later, still in the framework of \emph{finite} von Neumann algebras, it was studied in detail by P. Jolissaint in \cite{Jol}, where he proved in particular that in fact it does not depend on the choice of the faithful normal tracial state.

It is fair to say that the interest in extending the notion of the Haagerup property to infinite von Neumann algebras (or even just considering non-tracial states when building $L^2$-approximations on finite von Neumann algebras) was rather limited for several years due to the fact that in general there is no hope to characterise the Haagerup property for a non-discrete group via its von Neumann algebra. In the same spirit as the relation between  amenability of $G$ and injectivity of $VN(G)$ breaks down for non-discrete groups. This changed drastically in the last few years, with the advent of a study of the Haagerup property for discrete and locally compact \emph{quantum} groups. It was initiated by the study of the dual of the quantum free orthogonal and free unitary group by M. Brannan (\cite{Brannan}) and then continued in several directions (see \cite{DawFimSkaWhi} and references therein). The key factor lies in the fact that the natural (Haar) states on the von Neumann algebras of discrete quantum groups need not be tracial, so to characterise the Haagerup property for a discrete quantum group via its von Neumann algebra, one first needs to develop an understanding of the von Neumann algebraic Haagerup approximation property in absence of a tracial state.

The authors of this note approached this question from two directions: MC and AS in \cite{CasSka} studied a direct variant of the original definition, based on the existence of the approximating completely positive maps which behave well with respect to a given faithful normal weight and whose $L^2$-GNS implementations are compact, whereas RO and RT in \cite{OkaTom} worked in the framework of standard forms, partly motivated by the approach to injectivity due to A. Torpe (\cite{Torpe}). After the original versions of these preprints were circulated, the current authors observed that in fact the two approaches lead to equivalent notions (although this only becomes clear after the respective theories are fully developed). In this note we summarise the main results regarding the Haagerup property for arbitrary von Neumann algebras, referring to \cite{CasSka} and \cite{OkaTom} for the proofs of the presented statements.

We assume throughout that the von Neumann algebras we consider have separable preduals and follow the notation of \cite{TakesakiII}.

\section{Main results}

Let $\mlg$ be a von Neumann algebra and let $\varphi$ be a normal semifinite faithful weight on $\mlg$. The GNS Hilbert space arising as the completion of the set  $\nphi:=\{x\in \mlg:\varphi(x^*x) <\infty\}$ with respect to the norm associated to the scalar product $\la x, y \ra_{\varphi} := \varphi(y^*x)$ will be denoted by $L^2(\mlg,\varphi)$, and the GNS embedding map $\nphi \to L^2(\mlg,\varphi)$ by $\Lambda_{\varphi}$.
If $\Phi:\mlg \to \mlg$ is a normal, completely positive map which is \emph{$\varphi$-non-increasing}, i.e.\ it satisfies the inequality $\varphi(\Phi(x)) \leq \varphi(x)$ for all $x \in \mlg_+$, then $\Phi$ induces a bounded map on  $L^2(\mlg,\varphi)$ via the continuous linear extension of the prescription
$ \Lambda_{\varphi}(\nphi) \rightarrow  \Lambda_{\varphi}(\nphi): \Lambda_{\varphi}(x) \mapsto \Lambda_{\varphi}(\Phi(x))$. The following definition is considered in \cite{CasSka}.

\vspace*{0.2cm}

\begin{e-definition} \label{phiHAP}
Let $\mlg$ be a von Neumann algebra and let $\varphi$ be a normal semifinite faithful weight on $\mlg$. We say that $\mlg$  has the Haagerup approximation property with respect to $\varphi$ if there exists a sequence of completely positive, $\varphi$-non-increasing maps $(\Phi_n)_{n=1}^{\infty}$ on $\mlg$ such that the associated induced maps on the Hilbert space $L^2(\mlg, \varphi)$ are compact and converge to $I_{L^2(\mlg, \varphi)}$ strongly.
\end{e-definition}

\vspace*{0.2cm}





Recall that if $\mlg$ is a von Neumann algebra, then a quadruple $(\mlg, \Hil, J, \Pcone)$ is said to be a \emph{standard form} of $\mlg$ if $\Hil$ is a Hilbert space on which $\mlg$ is faithfully and non-degenerately represented, $J$ is a conjugate-linear isometry of $\Hil$  with $J^2 = 1$ and $\Pcone \subset \Hil$ a self-dual closed convex cone, such that $JMJ=M'$, $J$  fixes $\Pcone$ pointwise and $xJxJ$  fixes $\Pcone$ globally for each $x \in \mlg$. The existence (and uniqueness up to unitary equivalence) of standard forms was established in \cite{Haastf}.

If $(\mlg, \Hil, J, \Pcone)$ is a standard form of $\mlg$ and $T$ is a bounded linear operator acting on $\Hil$ then we say that $T$ is \emph{positive} if it preserves the cone $\Pcone$. Similarly, using the fact that if $n\in \bn$ then the standard form of the von Neumann algebra $M_n(\mlg)$ can be realised in the form $(M_n(\mlg), M_n(\Hil), J_n, \Pcone^{(n)})$, where $\Pcone^{(n)}$ is a certain cone in the Hilbert space $M_n (\Hil)$ (identified as the tensor product of $\Hil$ and $M_n$ equipped with the Hilbert-Schmidt norm), we can also say that $T\in B(\Hil)$ is  \emph{completely positive} if for all $n \in \bn$ the natural matrix lifting $T^{(n)}:M_n (\Hil) \to M_n (\Hil)$ preserves the cone $\Pcone^{(n)}$ (see \cite{MT} and \cite{SW}). The following notion of the Haagerup property is considered in \cite{OkaTom}.

\vspace*{0.2cm}

\begin{e-definition} \label{stfHAP}
Let $\mlg$ be a von Neumann algebra with a standard form $(\mlg, \Hil, J, \Pcone)$. We say that $\mlg$  has the standard form Haagerup approximation property if there exists a sequence of contractive completely positive compact operators $(T_n)_{n=1}^{\infty}$ on $\Hil$ such that $T_n \stackrel{n \to \infty}{\longrightarrow} I_{\Hil}$ strongly.
\end{e-definition}

\vspace*{0.2cm}

\begin{theorem} \label{equiv}
Let $\mlg$ be a von Neumann algebra with   separable predual. The following conditions are equivalent:
\begin{rlist}
\item $\mlg$  has the standard form Haagerup approximation property;
\item $\mlg$  has the Haagerup approximation property with respect to some  normal semifinite faithful weight $\varphi$;
\item $\mlg$  has the Haagerup approximation property with respect to every normal semifinite faithful weight  $\varphi$.
\end{rlist}
If $\mlg$ is a finite von Neumann algebra they are also equivalent to the Haagerup approximation property as defined and studied in \cite{Cho} and \cite{Jol}.
\end{theorem}

\vspace*{0.2cm}

In view of the above theorem and the comments before we  will simply say that a von Neumann algebra $\mlg$ with  separable predual has the \emph{Haagerup property} if it satisfies the equivalent conditions above. We will postpone the discussion of the proof of Theorem \ref{equiv} to the last part of the article.

\begin{theorem}
The following classes of (possibly infinite) von Neumann algebras have the Haagerup property:
\begin{itemize}
\item injective von Neumann algebras, in particular $B(\Hil$) (\cite[Corollary 2.10]{OkaTom});
\item von Neumann algebras of discrete quantum groups with the Haagerup property (\cite[Theorem 6.4]{DawFimSkaWhi}); in particular the von Neumann algebras of the duals of the quantum free orthogonal and unitary groups (\cite{ComFreYam}).
\end{itemize}
\end{theorem}

\vspace*{0.2cm}

The Haagerup property is also stable under several natural von Neumann algebraic constructions, which we list in the next theorem.

\vspace*{0.2cm}

\begin{theorem}
The Haagerup property is stable under the following constructions:
\begin{itemize}
\item tensor products: $\mlg_1 \overline{\ot} \mlg_2$ has the Haagerup property if and only if both $\mlg_1$ and $ \mlg_2$ have the Haagerup property \cite[Theorem 3.7]{OkaTom}, \cite[Lemma 3.5, Proposition 5.8]{CasSka};
\item direct sums: $\bigoplus_{n \in \bn} \mlg_n$ has the Haagerup property if and only if each $\mlg_n$ has the Haagerup property \cite[Theore 3.14]{OkaTom};
\item passing to the commutant: $\mlg$ has the Haagerup property if and only if $\mlg'$ does \cite[Theorem 3.12]{OkaTom};
\item passing to a corner: if $p \in \mlg$ is a projection and $\mlg$ has the Haagerup property then $p\mlg p$ has the Haagerup property \cite[Proposition 3.5]{OkaTom}, \cite[Proposition 5.9]{CasSka}. Moreover if  $(p_n)_{n=1}^{\infty}$ is an increasing sequence of projections in $\mlg$ strongly convergent to $1_{\mlg}$ and each $p_n\mlg p_n$ has the Haagerup property, then $\mlg$ has the Haagerup property \cite[Proposition 3.5]{OkaTom};
\item taking conditional expectations: if $\mlg$ has the Haagerup property and $\nlg$ is a von Neumann algebra which is the image of a conditional expectation $E:\mlg \to \nlg$ (not necessarily normal!) then $\nlg$ has the Haagerup property \cite[Theorem 5.9]{OkaTom}.
\end{itemize}
\end{theorem}

\vspace*{0.2cm}

Crucially, the Haagerup property behaves well under  crossed products, especially by amenable groups.

\vspace*{0.2cm}

\begin{theorem} \label{crossed}
Let $G$ be a locally compact group and let $\alpha:G \to \textup{Aut} (\mlg)$ be a continuous action of $G$ on the von Neumann algebra $\mlg$ (this means that $\alpha$ is a homomorphism and the map
$G \ni g \mapsto \alpha_g(x) \in \mlg$ is strongly continuous for each $x \in \mlg$). Then the following hold:
\begin{rlist}
\item if $\mlg \rtimes_{\alpha} G$  has the  Haagerup property, then $\mlg$ has the Haagerup property \cite[Corollary 5.14]{OkaTom}, \cite[Section 5]{CasSka};
\item if $G$ is amenable and $\mlg$ has the Haagerup property, then $\mlg \rtimes_{\alpha} G$  has the  Haagerup property \cite[Corollary 5.15]{OkaTom}, \cite[Section 6]{CasSka}.
\end{rlist}
\end{theorem}

\vspace*{0.2cm}

The latter result can also be generalised to crossed products by the actions of locally compact amenable quantum groups with amenable duals, see the same references.

\vspace*{0.2cm}

Recall that the \emph{core}  of a von Neumann algebra $\mlg$ is the crossed product of $\mlg$ by the action of the modular automorphism of a normal semifinite faithful weight on $\mlg$ (it is known to be independent of the choice of such a weight). Theorem \ref{crossed} implies directly the following corollary.

\vspace*{0.2cm}

\begin{corollary} \label{core}
A von Neumann algebra has the Haagerup property if and only if its core does.
\end{corollary}

\vspace*{0.2cm}

Corollary \ref{core} is established separately for the Haagerup approximation property with respect to a normal semifinite faithful weight \cite{CasSka} and for the standard form Haagerup approximation property \cite{OkaTom} (see Definitions \ref{phiHAP} and \ref{stfHAP}). This forms the key to the proof of Theorem \ref{equiv}. Following \cite{CasSka} the passage to the semifinite core allows one to prove the implication (ii)$\Longrightarrow$(iii) in Theorem \ref{equiv}. On the other hand, since the core of any von Neumann $\mlg$ can be decomposed as a tensor product of a finite von Neumann algebra $\nlg$ and a I$_\infty$-factor, this shows that for both the Definitions \ref{phiHAP} and \ref{stfHAP} the Haagerup property of $\mlg$ is equivalent to the Haagerup property of $\nlg$ which in turn agrees with the existing notion of the Haagerup approximation property in the case of finite von Neumann algebras.


Finally, another alternative definition of the Haagerup property is given in  \cite[Definition 3.7]{OkaTom2} case $\alpha=0$, replacing strong convergence of the $L^2$-maps in Definition \ref{phiHAP} by $\sigma$-weak convergence of the maps $\Phi_n$ on $\mlg$.  By  \cite[Theorem A]{OkaTom2} this definition is equivalent to the standard form Haagerup approximation property.

\end{document}